\documentclass{amsart}
\usepackage{amssymb, amsmath, amsthm, graphics, comment, xspace, enumerate}
\newtheorem{thm}{Theorem}[section]

\newtheorem{prop}[thm]{Proposition}
\theoremstyle{definition}
\newtheorem{defn}[thm]{Definition}

\theoremstyle{remark}
\newtheorem{rem}[thm]{Remark}
\numberwithin{equation}{section}
\begin{document}
\title[Degenerate $C$-ultradistribution semigroups in lcs]{Degenerate $C$-ultradistribution semigroups in locally convex spaces}
\author{Marko Kosti\' c}
\address{Faculty of Technical Sciences,
University of Novi Sad,
Trg D. Obradovi\' ca 6, 21125 Novi Sad, Serbia}
\email{marco.s@verat.net}

\author{Stevan Pilipovi\' c}
\address{Department for Mathematics and Informatics,
University of Novi Sad,
Trg D. Obradovi\' ca 4, 21000 Novi Sad, Serbia}
\email{pilipovic@dmi.uns.ac.rs}

\author{Daniel Velinov}
\address{Department for Mathematics, Faculty of Civil Engineering, Ss. Cyril and Methodius University, Skopje,
Partizanski Odredi
24, P.O. box 560, 1000 Skopje, Macedonia}
\email{velinovd@gf.ukim.edu.mk}

{\renewcommand{\thefootnote}{} \footnote{
\\ \text{  }  \ \    {\it Key words and phrases.} Degenerate $C$-ultradistribution semigroups, multivalued linear operators, locally convex spaces.
\\  \text{  }  \ \ This research is partially supported by grant 174024 of Ministry
of Science and Technological Development, Republic of Serbia.}}

\begin{abstract}
The main subject in this paper are degenerate $C$-ultradistribution semigroups in barreled
sequentially complete locally convex spaces.
Here, the regularizing operator $C$ is not necessarily injective and the infinitesimal generator is multivalued linear operator. We also consider exponential degenerate $C$-ultradistribution semigroups.
\end{abstract}
\maketitle

\section{Introduction and Preliminaries}
This is an expository paper. We collect known results and the results which simply follows from the known ones. Because of that proofs are not given. In \cite{C-ultra} are introduced and systematically analyzed
the classes of $C$-distribution and $C$-ultradistribution semigroups in locally convex spaces (cf. \cite{b42}-\cite{cizi}, \cite{fat1},
\cite{ki90}, \cite{k92}-\cite{knjigaho}, \cite{ku112}-\cite{li121},
\cite{me152}, \cite{1964}-\cite{w241}  and references cited therein for the current state of theory).
The recent paper \cite{degdis} motivate us to continue the study on generalized degenerate $C$-regularized semigroups in locally convex spaces in the case of ultradistribution semigroups.
The main aim of this paper is to investigate the degenerate $C$-ultradistribution semigroups in the setting of barreled
sequentially complete locally convex spaces. We refer to \cite{carol}, \cite{faviniyagi}, \cite{FKP},
\cite{me152} and \cite{svir-fedorov} for further information about well-posedness of abstract degenerate differential equations of first order. Here, we consider multivalued linear operators as
infinitesimal generators of a degenerate $C$-ultradistribution semigroups (cf. \cite{baskakov-chern}, \cite{ki90}, \cite{ku112}, \cite{isna-maiz}). 
%Like in \cite{degdis}, we do not use any decomposition of the state space
%$E.$
The organization of the paper is as follows. In Section 1 are exposed the basic facts about vector-valued ultradistributions. Our main results are contained in Section 2, in which we analyze various themes concerning degenerate $C$-ultradistribution semigroups in locally convex spaces and further generalize some of our recent results from \cite{C-ultra} and \cite{degdis}.

\subsection{Notation}
We use the standard notation throughout the paper.
Unless specified otherwise,
we assume
that $E$ is a Hausdorff sequentially complete
locally convex space over the field of complex numbers, SCLCS for short.  For the sake of brevity
and better exposition, our standing assumption henceforth will be that
the state space $ E$ is barreled.
By
$L(E)$ we denote the space consisting of all continuous linear mappings from $E$ into
$E$ and by the symbol $\circledast_{E}$ (usually we will denote $\circledast$ if there is no risk for confusion) denotes the fundamental system of seminorms which defines the topology of $E.$
Let $X$ be also an SCLCS, let ${\mathcal B}$ be the family of bounded subsets\index{bounded subset} of $E,$ and
let $p_{B}(T):=\sup_{x\in B}p(Tx),$ $p\in \circledast_{X},$ $B\in
{\mathcal B},$ $T\in L(E,X).$ Then $p_{B}(\cdot)$ is a seminorm\index{seminorm} on
$L(E,X)$ and the system $(p_{B})_{(p,B)\in \circledast_{X} \times
{\mathcal B}}$ induces the Hausdorff locally convex topology on
$L(E,X).$
The Hausdorff locally convex topology on $E^{\ast},$ the dual space\index{dual space} of $E,$
defines the
system $(|\cdot|_{B})_{B\in {\mathcal B}}$ of seminorms on
$E^{\ast},$ where $|x^{\ast}|_{B}:=\sup_{x\in
B}|\langle x^{\ast}, x \rangle |,$ $x^{\ast} \in E^{\ast},$ $B\in
{\mathcal B}.$ The bidual of $E$ is denoted by $E^{\ast \ast}.$ The polars of nonempty sets $M\subseteq E$ and $N\subseteq E^*$ are defined as follows
$M^{\circ}:=\{y\in E^*:|y(x)|\leq 1\text{ for all } x\in M\}$ and
$N^{\circ}:=\{x\in E:\;|y(x)|\leq 1\text{ for all } y\in N\}.$
If $A$ is a linear operator
acting on $E,$
then the domain, kernel space and range of $A$ will be denoted by
$D(A),$ $N(A)$ and $R(A),$
respectively. Since no confusion
seems likely, we will identify $A$ with its graph.
Since we have assumed that the state space $ E$ is barreled, the spaces $L(E)$ and $E^{\ast}$ are sequentially
complete (\cite{meise}) and any strongly continuous operator family
$(S(t))_{t\in [0,\tau)}\subseteq L(E),$ where $0 <\tau \leq \infty,$  is locally equicontinuous. The reader may consult \cite{x263} and \cite{FKP}
for further information on the Laplace transform of functions with values in SCLCS's; cf. \cite{a43} for the Banach space case.\\
We assume that $(M_p)$ is a sequence of positive real numbers
such that $M_0=1$ and the following conditions hold:\vspace{0.1cm} \newline \noindent
(M.1): $M_p^2\leq M_{p+1} M_{p-1},\;\;p\in\mathbb{N},
$\\
(M.2): $M_p\leq AH^p\sup_{0\leq i\leq p}M_iM_{p-i},\;\;p\in\mathbb{N},\mbox{ for some }A,\ H>1,
$\\
(M.3)': $\sum_{p=1}^{\infty}\frac{M_{p-1}}{M_p}<\infty.
$\\
Every employment of the condition\vspace{0.1cm} \newline \noindent
(M.3): $\sup_{p\in\mathbb{N}}\sum_{q=p+1}^{\infty}\frac{M_{q-1}M_{p+1}}{pM_pM_q}<\infty,
$\\
which is a slightly stronger than (M.3)', will be explicitly emphasized.

Let $s>1$.
Then the Gevrey sequence $(p!^s)$ satisfies the above conditions.
The associated
function of sequence $(M_p)$ is defined by
$M(\rho):=\sup_{p\in\mathbb{N}}\ln\frac{\rho^p}{M_p}$,
$\rho>0$; $M(0):=0,$ $M(\lambda):=M(|\lambda|),$
$\lambda\in\mathbb{C} \setminus [0,\infty).$

Let us recall that the spaces of Beurling,
respectively, Roumieu ultradifferentiable functions are
defined by $\mathcal{D}^{(M_p)}:=\mathcal{D}^{(M_p)}(\mathbb{R})
:=\text{indlim}_{K\Subset\Subset\mathbb{R}}\mathcal{D}^{(M_p)}_K$,
respectively,
$\mathcal{D}^{\{M_p\}}:=\mathcal{D}^{\{M_p\}}(\mathbb{R})
:=\text{indlim}_{K\Subset\Subset\mathbb{R}}\mathcal{D}^{\{M_p\}}_K$, (where $K$ goes through all compact sets in ${\mathbb R}$
where
$\mathcal{D}^{(M_p)}_K:=\text{projlim}_{h\to\infty}\mathcal{D}^{M_p,h}_K$,
respectively, $\mathcal{D}^{\{M_p\}}_K:=\text{indlim}_{h\to 0}\mathcal{D}^{M_p,h}_K$,
\begin{align*}
\mathcal{D}^{M_p,h}_K:=\bigl\{\phi\in C^{\infty}(\mathbb{R}): \text{supp}(\phi) \subseteq K,\;\|\phi\|_{M_p,h,K}<\infty\bigr\},
\end{align*}

\begin{align*}
\|\phi\|_{M_p,h,K}:=\sup\Biggl\{\frac{h^p\bigl|\phi^{(p)}(t)\bigr|}{M_p} : t\in K,\;p\in\mathbb{N}_0\Biggr\}.
\end{align*}
Spaces of tempered ultradistributions
are defined as strong dual of corresponding test spaces:
$$
{\mathcal S}^{(M_p)}({\mathbb R}):=\mbox{proj}\lim_{k\rightarrow
\infty}{\mathcal S}^{M_p,k}({\mathbb R}), \mbox{ resp., }
{\mathcal S}^{\{M_p\}}({\mathbb R}):=\mbox{ind}
\lim_{k\rightarrow 0}{\mathcal S}^{M_p,k}({\mathbb R}),
$$
where
$${\mathcal S}^{M_p,k}({\mathbb R}):=\{\phi\in
C^\infty({\mathbb R}) : ||\phi||_{M_p,k}<\infty\},\ k>0,$$
$$||\phi||_{M_p,k}:=\mbox{ sup}\{\frac{k^{\alpha+\beta}}{M_\alpha
M_\beta}(1+|t|^2)^{\beta/2}|\phi^{(\alpha)}(t)| :  t\in {\mathbb
R},\ \alpha,\ \beta \in {{\mathbb N}_0} \}.$$
Henceforth the asterisk $*$ stands for both cases. Let $\emptyset \neq \Omega \subseteq {\mathbb R}.$
The spaces
$\mathcal{D}'^*(E):=L(\mathcal{D}^*, E)$, $\mathcal{D}^{*}_{\Omega}$, $\mathcal{D}^{\ast}_0$, $\mathcal{E}'^{*}_{\Omega}$, $\mathcal{E}'^{*}_{0}$, $\mathcal{D}'^{*}_{\Omega}(E)$, $\mathcal{D}'^{*}_{0}(E)$ and $\mathcal{S}'^{\ast}_0(E)$ are defined as in distribution case.
%For every $t\in {\mathbb R},$ we define the scalar-valued ultradistribution
%$\delta_{t}\in \mathcal{D}'^{*}$ by $\delta_{t}(\varphi):=\varphi(t),$ $\varphi \in {\mathcal D^{*}}$.
%The multiplication by a function $a\in {\mathcal E}^{\ast}(\Omega)$,
%convolution of scalar valued ultradistributions (ultradifferentiable
%functions), and the notion of a regularizing sequence in $\mathcal{D}^*,$ are defined as in the case of
%distributions;
We know that there exists
a regularizing sequence in $\mathcal{D}^*$.
Regularizing
sequence in $\mathcal{D}^{\ast}$ we mean is a sequence $(\rho_n)_{n\in {\mathbb N}}$ in
$\mathcal{D}_0^{\ast}$ for which there exists a function $\rho\in\mathcal{D}^{\ast}$ such that $\int_{-\infty}^{\infty}\rho
(t)\,dt=1,$ supp$(\rho)\subseteq [0,1]$ and $\rho_n(t)=n\rho(nt)$,
$t\in\mathbb{R},$ $n\in {\mathbb N}.$
We define the convolution products $\varphi*\psi$
and $\varphi*_0\psi$ by
$$
\varphi*\psi(t):=\int\limits_{-\infty}^{\infty}\varphi(t-s)\psi(s)\,ds\mbox{ and }
\varphi*_0
\psi(t):=\int\limits^t_0\varphi(t-s)\psi(s)\,ds,\;t\;\in\mathbb{R},
$$ for $\varphi$, $\psi:\mathbb{R}\to\mathbb{C}$
locally integrable functions.
Notice that $\varphi*\psi=\varphi*_0\psi$, provided that supp$(\varphi)$ and supp$(\psi)$ are subsets of $[0,\infty).$
Given $\varphi\in\mathcal{D}^{*}$ and $f\in\mathcal{D}'^{\ast}$, or $\varphi\in\mathcal{E}^{*}$ and $f\in\mathcal{E}'^{*}$,
we define the convolution $f*\varphi$ by $(f*\varphi)(t):=f(\varphi(t-\cdot))$, $t\in\mathbb{R}$.
For $f\in\mathcal{D}'^{\ast}$, or for $f\in\mathcal{E}'^{\ast}$,
define $\check{f}$ by $\check{f}(\varphi):=f(\varphi (-\cdot))$, $\varphi\in\mathcal{D}^{\ast}$ ($\varphi\in\mathcal{E}^{*}$).
The convolution of two ultradistributions $f$, $g\in\mathcal{D}'^{\ast}$, denoted by $f*g$,
is defined by $(f*g)(\varphi):=g(\check{f}*\varphi)$, $\varphi\in\mathcal{D}^{\ast}$.\\
%If one of them belongs to ${\mathcal E}'({\mathbb R})$, then we know that $f*g\in\mathcal{D}'$ and supp$(f*g)\subseteq$supp$ (f)+$supp$(g)$.
%If $\varphi\in\mathcal{D}^{\ast}$ ($T\in {\mathcal E}'^{\ast}$) and $G\in\mathcal{D}'^{\ast}(E)$, then
%$\varphi \ast G \in {\mathcal E}^{\prime \ast}(E)$ and $T \ast G \in  {\mathcal D}^{\prime \ast}(E)$ (cf. \cite[p. 685]{k82}, and
%\cite[Definition 3.9]{k82} for the notion of space ${\mathcal E}^{\prime \ast}(E)$).

%We refer the reader to \cite{C-ultra} for some characterizations of vector-valued (ultra)distributions supported by a point.
%If the space $E$ satisfies the property that any vector-valued distribution $G\in\mathcal{D}'(E)$ with supp$(G)\subseteq\{0\}$ can be represented as a finite sum of
%vector-valued distributions of form $\delta^{(i)} \otimes x_{i}$, then we
%say that $E$ is admissible.

%\section[Multivalued linear operators]{Multivalued linear operators}

%In this section, we present some definitions and properties of multivalued linear operators that will be necessary for our further work (cf. the monographs \cite{cross} by R. Cross and \cite{faviniyagi} by A. Favini-A. Yagi for more details on the subject).
%The underlying SCLCS will be denoted
%by $X$ and $Y;$ in the third section, we will coming back to our standing notation.
\indent We recall the definition of a multivalued map (multimap) given in our recent paper \cite{degdis} (cf. \cite{cross} by R. Cross, \cite{faviniyagi} by A. Favini-A. Yagi).
A multivalued map (multimap) ${\mathcal A} : X \rightarrow P(Y)$ is said to be a multivalued
linear operator (MLO) iff the following holds:
\begin{itemize}
\item[(i)] $D({\mathcal A}) := \{x \in X : {\mathcal A}x \neq \emptyset\}$ is a subspace of $X$;
\item[(ii)] ${\mathcal A}x +{\mathcal A}y \subseteq {\mathcal A}(x + y),$ $x,\ y \in D({\mathcal A})$
and $\lambda {\mathcal A}x \subseteq {\mathcal A}(\lambda x),$ $\lambda \in {\mathbb C},$ $x \in D({\mathcal A}).$
\end{itemize}
If $X=Y,$ then it is also said that ${\mathcal A}$ is an MLO in $X.$
%An almost immediate consequence of the definition is that,
%for every $x,\ y\in D({\mathcal A})$ and for every $\lambda,\ \eta \in {\mathbb C}$ with $|\lambda| + |\eta| \neq 0,$ we
%have $\lambda {\mathcal A}x + \eta {\mathcal A}y = {\mathcal A}(\lambda x + \eta y).$ If ${\mathcal A}$ is an MLO, then ${\mathcal A}0$ is a linear manifold in $Y$
%and ${\mathcal A}x = f + {\mathcal A}0$ for any $x \in D({\mathcal A})$ and $f \in {\mathcal A}x.$ Set $R({\mathcal A}):=\{{\mathcal A}x :  x\in D({\mathcal A})\}.$
%The set ${\mathcal A}^{-1}0 = \{x \in D({\mathcal A}) : 0 \in {\mathcal A}x\}$ is called the kernel\index{multivalued linear operator!kernel}
%of ${\mathcal A}$ and it is denoted by $N({\mathcal A}).$
The inverse ${\mathcal A}^{-1}$ of an MLO is defined by
$D({\mathcal A}^{-1}) := R({\mathcal A})$ and ${\mathcal A}^{-1} y := \{x \in D({\mathcal A}) : y \in {\mathcal A}x\}$.\index{multivalued linear operator!inverse}
It is easily seen that ${\mathcal A}^{-1}$ is an MLO in $X,$ as well as that $N({\mathcal A}^{-1}) = {\mathcal A}0$
and $({\mathcal A}^{-1})^{-1}={\mathcal A}.$ If $N({\mathcal A}) = \{0\},$ i.e., if ${\mathcal A}^{-1}$ is
single-valued, then ${\mathcal A}$ is said to be injective.

%For any mapping ${\mathcal A}: X \rightarrow P(Y)$ we define $\check{{\mathcal A}}:=\{(x,y) : x\in D({\mathcal A}),\ y\in {\mathcal A}x\}.$ Then ${\mathcal A}$ is an MLO iff $\check{{\mathcal A}}$ is a linear relation in $X\times Y,$ ($(x,\lambda y_1)+(x,\lambda y_2)=(x,\lambda y_1+\lambda y_2)$, for $x\in X$ and $y\in Y$) i.e., iff $\check{{\mathcal A}}$ is a subspace of $X \times Y.$ Since no confusion
%seems likely, we will sometimes identify ${\mathcal A}$ with its graph. \index{linear relation}

If ${\mathcal A},\ {\mathcal B} : X \rightarrow P(Y)$ are two MLOs, then we define its sum ${\mathcal A}+{\mathcal B}$ by $D({\mathcal A}+{\mathcal B}) := D({\mathcal A})\cap D({\mathcal B})$ and $({\mathcal A}+{\mathcal B})x := {\mathcal A}x +{\mathcal B}x,$ $x\in D({\mathcal A}+{\mathcal B}).$
It can be simply checked that ${\mathcal A}+{\mathcal B}$ is likewise an MLO.\index{multivalued linear operator!sum}

Let ${\mathcal A} : X \rightarrow P(Y)$ and ${\mathcal B} : Y\rightarrow P(Z)$ be two MLOs, where $Z$ is an SCLCS. The product of ${\mathcal A}$
and ${\mathcal B}$ is defined by $D({\mathcal B}{\mathcal A}) :=\{x \in D({\mathcal A}) : D({\mathcal B})\cap {\mathcal A}x \neq \emptyset\}$ and\index{multivalued linear operator!product}
${\mathcal B}{\mathcal A}x:=
{\mathcal B}(D({\mathcal B})\cap {\mathcal A}x).$ Then ${\mathcal B}{\mathcal A} : X\rightarrow P(Z)$ is an MLO and
$({\mathcal B}{\mathcal A})^{-1} = {\mathcal A}^{-1}{\mathcal B}^{-1}.$ The scalar multiplication of an MLO ${\mathcal A} : X\rightarrow P(Y)$ with the number $z\in {\mathbb C},$ $z{\mathcal A}$ for short, is defined by
$D(z{\mathcal A}):=D({\mathcal A})$ and $(z{\mathcal A})(x):=z{\mathcal A}x,$ $x\in D({\mathcal A}).$ It is clear that $z{\mathcal A}  : X\rightarrow P(Y)$ is an MLO and $(\omega z){\mathcal A}=\omega(z{\mathcal A})=z(\omega {\mathcal A}),$ $z,\ \omega \in {\mathbb C}.$

The integer powers of an MLO ${\mathcal A} :  X\rightarrow P(X)$ is defined recursively as follows: ${\mathcal A}^{0}=:I;$ if ${\mathcal A}^{n-1}$ is defined, set $
D({\mathcal A}^{n}) := \bigl\{x \in  D({\mathcal A}^{n-1}) : D({\mathcal A}) \cap {\mathcal A}^{n-1}x \neq \emptyset \bigr\},
$
and
$
{\mathcal A}^{n}x := \bigl({\mathcal A}{\mathcal A}^{n-1}\bigr)x =\bigcup_{y\in  D({\mathcal A}) \cap {\mathcal A}^{n-1}x}{\mathcal A}y,\quad x\in D( {\mathcal A}^{n}).
$
It is well known that $({\mathcal A}^{n})^{-1} = ({\mathcal A}^{n-1})^{-1}{\mathcal A}^{-1} = ({\mathcal A}^{-1})^{n}=:{\mathcal A}^{-n},$ $n \in {\mathbb N}$
and $D((\lambda-{\mathcal A})^{n})=D({\mathcal A}^{n}),$ $n \in {\mathbb N}_{0},$ $\lambda \in {\mathbb C}.$ Moreover,
if ${\mathcal A}$ is single-valued, then the above definitions are consistent with the usual definition of powers of ${\mathcal A}.$

%If ${\mathcal A} : X\rightarrow P(Y)$ and ${\mathcal B} : X\rightarrow P(Y)$ are two MLOs, then we write ${\mathcal A} \subseteq {\mathcal B}$ iff $D({\mathcal A}) \subseteq D({\mathcal B})$ and ${\mathcal A}x \subseteq {\mathcal B}x$
%for all $x\in D({\mathcal A}).$ Assume now that
%a linear single-valued operator $S : D(S) \subseteq X \rightarrow Y$ has domain $D(S) = D({\mathcal A})$ and $S \subseteq {\mathcal A},$ where ${\mathcal A} : X\rightarrow P(Y)$
%is an MLO. Then $S$ is called a\index{multivalued linear operator!section}
%section of ${\mathcal A};$ if this is the case, we have ${\mathcal A}x = Sx + {\mathcal A}0,$ $x \in D({\mathcal A})$ and
%$R({\mathcal A}) = R(S) + {\mathcal A}0.$
%
%We say that an MLO operator  ${\mathcal A} : X\rightarrow P(Y)$ is closed if for any
%nets $(x_{\tau})$ in $D({\mathcal A})$ and $(y_{\tau})$ in $Y$ such that $y_{\tau}\in {\mathcal A}x_{\tau}$ for all $\tau\in I$ we have that $\lim_{\tau \rightarrow \infty}x_{\tau}=x$ and
%$\lim_{\tau \rightarrow \infty}y_{\tau}=y$ imply
%$x\in D({\mathcal A})$ and $y\in {\mathcal A}x.$\index{multivalued linear operator!closed}

If ${\mathcal A} : X\rightarrow P(Y)$ is an MLO, then we define the adjoint ${\mathcal A}^{\ast}: Y^{\ast}\rightarrow P(X^{\ast})$\index{multivalued linear operator!adjoint}
of ${\mathcal A}$ by its graph
$$
{\mathcal A}^{\ast}:=\Bigl\{ \bigl( y^{\ast},x^{\ast}\bigr)  \in Y^{\ast} \times X^{\ast} :  \bigl\langle y^{\ast},y \bigr\rangle =\bigl \langle x^{\ast}, x\bigr \rangle \mbox{ for all pairs }(x,y)\in {\mathcal A} \Bigr\}.
$$
%It is simply verified that ${\mathcal A}^{\ast}$
%is a closed MLO, and that $ \langle y^{\ast},y \rangle =0$ whenever $y^{\ast}\in D({\mathcal A}^{\ast})$ and $y\in {\mathcal A}0.$

%Concerning the integration of functions with values in SCLCS, we follow the approach of C. Martinez and M. Sanz \cite[pp. 99-102]{martinez}. Denote by $\Omega$ a locally compact\index{space!locally compact} and separable metric space\index{space!separable metric}
%and  by $\mu$ a locally finite
%Borel measure\index{measure!locally finite Borel} defined on $\Omega.$ Then the following fundamental lemma holds:
%
%\begin{lem}\label{integracija-tricky}
%Suppose that ${\mathcal A} : X\rightarrow P(Y)$ is a closed \emph{MLO}. Let $f : \Omega \rightarrow X$ and $g : \Omega \rightarrow Y$ be $\mu$-integrable, and let $g(x)\in {\mathcal A}f(x),$ $x\in \Omega.$ Then $\int_{\Omega}f\, d\mu \in D({\mathcal A})$ and $\int_{\Omega}g\, d\mu\in {\mathcal A}\int_{\Omega}f\, d\mu.$
%\end{lem}

In \cite{FKP}, we have recently considered the  $C$-resolvent sets of MLOs in locally convex spaces
(where $C\in L(X)$ is injective, $C{\mathcal A}\subseteq {\mathcal A}C$).
The
$C$-resolvent set of an MLO ${\mathcal A}$ in $X,$ $\rho_{C}({\mathcal A})$ for short, is defined as the union of those complex numbers
$\lambda \in {\mathbb C}$ for which $R(C)\subseteq R(\lambda-{\mathcal A})$ and
$(\lambda - {\mathcal A})^{-1}C$ is a single-valued bounded operator on $X.$
The operator $\lambda \mapsto (\lambda -{\mathcal A})^{-1}C$ is called the $C$-resolvent of ${\mathcal A}$ ($\lambda \in \rho_{C}({\mathcal A})$).
Here, we analyze the general situation in which the operator
$C\in L(X)$ is not necessarily injective. Then the operator $(\lambda - {\mathcal A})^{-1}C$ is no longer single-valued, which additionally hinders our considerations and work.

\section{Properties of the degenerate $C$-ultradistribution semigroups in locally convex spaces}

Throughout this section, we assume that $C\in L(E)$ is not necessarily injective operator. Since $E$ is barreled, the uniform boundedness principle \cite[p. 273]{meise} implies that each ${\mathcal G}\in {\mathcal D}'^{*}(L(E))$ is boundedly equicontinuous, i.e., that for every $p\in \circledast$ and for every bounded subset $B$ of ${\mathcal D}^{\ast}$, there exist $c>0$ and
$q\in \circledast$ such that
$
p({\mathcal G}(\varphi)x)\leq cq(x),\ \varphi \in B, \ x\in E.
$

We start this section by introducing the following definition.

\begin{defn}\label{cuds}
Let $\mathcal{G}\in\mathcal{D}_0'^{\ast}(L(E))$ satisfy $C\mathcal{G}=\mathcal{G}C.$
Then it is said that $\mathcal{G}$ is a pre-(C-UDS) of $\ast$-class iff the following holds:
\[\tag{C.S.1}
\mathcal{G}(\varphi*_0\psi)C=\mathcal{G}(\varphi)\mathcal{G}(\psi),\quad \varphi,\;\psi\in\mathcal{D}^{\ast}.
\]
If, additionally,
\[\tag{C.S.2}
\mathcal{N}(\mathcal{G}):=\bigcap_{\varphi\in\mathcal{D}^{\ast}_0}N(\mathcal{G}(\varphi))=\{0\},
\]
then $\mathcal{G}$ is called a $C$-ultradistribution semigroup of $\ast$-class, (C-UDS) of $\ast$-class in short.
A pre-(C-UDS) $\mathcal{G}$ is called dense iff
\[\tag{C.S.3}
\mathcal{R}(\mathcal{G}):=\bigcup\limits_{\varphi\in\mathcal{D}^{\ast}_0}R(\mathcal{G}(\varphi))
\text{ is dense in }E.
\
\]
\end{defn}

If $C=I,$ then we also write pre-(UDS), (UDS), instead of pre-(C-UDS), (C-UDS).

Suppose that $\mathcal{G}$ is a pre-(C-UDS) of $\ast$-class. Then
$\mathcal{G}(\varphi)\mathcal{G}(\psi)=\mathcal{G}(\psi)\mathcal{G}(\varphi)$ for all $\varphi,\,\psi\in\mathcal{D}^{\ast}$
and $\mathcal{N}(\mathcal{G})$ is a closed subspace of $E$.

The structural characterization of a pre-(C-UDS) $\mathcal{G}$ of $\ast$-class on its kernel space
$\mathcal{N}(\mathcal{G})$ is described in the following theorem (cf.  \cite[Proposition 3.1.1]{knjigah} and the proofs of \cite[Lemma 2.2]{ku112}, \cite[Proposition 3.5.4]{knjigah}).

\begin{thm}\label{delta-point}
Let $(M_{p})$ satisfy \emph{(M.3)}, let $\mathcal{G}$ be a pre-$($C-UDS$)$ of $\ast$-class, and let the space $\mathcal{N}(\mathcal{G})$ be barreled.
Then, with $N=\mathcal{N}(\mathcal{G})$ and $G_1$ being the restriction of $\mathcal{G}$ to $N$ $(G_1=\mathcal{G}_{|N})$,
we have:
There exists a unique set of operators $(T_{j})_{j\in {\mathbb N}_{0}}$ in $L(\mathcal{N}(\mathcal{G}))$  commuting with $C$ so that
$G_1=\sum_{j=0}^{\infty}\delta^{(j)}\otimes T_j$, $T_jC^j=(-1)^jT_0^{j+1}$, $j\in {\mathbb N}$ and the set $\{M_{j}T_{j}L^{j} : j\in {{\mathbb N}_{0}}\}$ is bounded in
$L(\mathcal{N}(\mathcal{G})),$ for some $L>0$ in the Beurling case, resp. for every $L>0$ in the Roumieu case.
\end{thm}

Let $\mathcal{G}\in\mathcal{D}_0'^{\ast}(L(E))$, and let $T\in\mathcal{E}_0'^{\ast}$,
i.e., $T$ is a scalar-valued ultradistribution of $\ast$-class with compact support contained in $[0,\infty)$.
Define
\[
G(T):=\Bigl\{(x,y) \in E\times E : \mathcal{G}(T*\varphi)x=\mathcal{G}(\varphi)y\;\mbox{ for all }\;\varphi\in\mathcal{D}^{\ast}_{0} \Bigr\}.
\]
Then  it can be easily seen that $G(T)$ is a closed MLO; furthermore, if $\mathcal{G}\in\mathcal{D}_0'^{\ast}(L(E))$ satisfy (C.S.2), then $G(T)$ is a closed linear operator.
Assuming that the regularizing operator $C$ is injective, definition of $G(T)$ can be equivalently introduced by replacing the set
$\mathcal{D}^{\ast}_{0}$ with the set $\mathcal{D}^{\ast}_{[0,\epsilon)}$ for any $\epsilon>0.$ In general case,
for every $\psi\in\mathcal{D}^{\ast}$, we have $\psi_+:=\psi\mathbf{1}_{[0,\infty)}\in\mathcal{E}_{0}'^{*}$, where
$\mathbf{1}_{[0,\infty)}$ stands for the characteristic function of $[0,\infty),$ so that the definition of $G(\psi_+)$ is clear.
We define the (infinitesimal) generator of a pre-(C-UDS) $\mathcal{G}$ by ${\mathcal A}:=G(-\delta')$ (cf. \cite{C-ultra} for more details about non-degenerate case, and \cite[Definition 3.4]{baskakov-chern} and \cite{ki90} for some other approaches
used in degenerate case). Then $\mathcal{N}(\mathcal{G}) \times \mathcal{N}(\mathcal{G}) \subseteq {\mathcal A}$ and $\mathcal{N}(\mathcal{G}) = {\mathcal A}0,$ which simply implies that
${\mathcal A}$ is single-valued iff (C.S.2)  holds. If this is the case, then we also have that the operator $C$ must be injective:
Suppose that $Cx=0$ for some $x\in E.$ By (C.S.1), we get that $\mathcal{G}(\varphi)\mathcal{G}(\psi)x=0,$ $\varphi,\;\psi\in\mathcal{D}.$ In particular, $\mathcal{G}(\psi)x\in {\mathcal N}({\mathcal G})=\{0\}$ so that $\mathcal{G}(\psi)x= 0,$ $\psi\in\mathcal{D}.$ Hence,  $x\in {\mathcal N}({\mathcal G})=\{0\}$ and therefore $x=0.$

Further on, if $\mathcal{G}$ is a pre-(C-UDS) of $\ast$-class, $T\in\mathcal{E}_{0}'^{*}$ and $\varphi\in\mathcal{D}^{\ast}$,
then ${\mathcal G}(\varphi)G(T)\subseteq G(T)\mathcal{G}(\varphi)$, $CG(T)\subseteq G(T)C$
and $\mathcal{R}(\mathcal{G})\subseteq D(G(T))$.
If $\mathcal{G}$ is a pre-(C-UDS) of $\ast$-class and $\varphi$, $\varphi$, $\psi\in\mathcal{D}^{\ast}$,
then the assumption $\varphi(t)=\psi(t)$, $t\geq 0$, implies $\mathcal{G}(\varphi)=\mathcal{G}(\psi)$.
As in the Banach space case, we can prove the following (cf. \cite[Proposition 3.1.3, Lemma 3.1.6]{knjigah}): Suppose that $\mathcal{G}$ is a pre-(C-UDS) of $\ast$-class. Then $(Cx,\mathcal{G}(\psi)x)\in G(\psi_+)$, $\psi\in\mathcal{D}^{\ast},$ $x\in E$  and ${\mathcal A}\subseteq C^{-1}{\mathcal A}C,$ while $ C^{-1}{\mathcal A}C={\mathcal A}$ provided that $C$ is injective.
The following two propositions holds in degenerate $C$-ultradistribution case (see \cite{degdis} for degenerate $C$-distribution case). Note that the reflexivity of the space $E$ implies that the spaces $E^{\ast}$ and $E^{\ast\ast}=E$ are both barreled and sequentially complete locally convex spaces.

\begin{prop}\label{isto}
Let ${\mathcal G}$ be a pre-(C-UDS) of $\ast$-class, $S$, $T\in\mathcal{E}'^{\ast}_0$, $\varphi\in\mathcal{D}^{\ast}_0$,
$\psi\in\mathcal{D}^{\ast}$ and $x\in E$.
Then we have:
\begin{itemize}
\item[(i)] $(\mathcal{G}(\varphi)x$, $\mathcal{G}(\overbrace{T*\cdots*T}^m*\varphi)x)\in G(T)^m$, $m\in\mathbb{N}$.
\item[(ii)] $G(S)G(T)\subseteq G(S*T)$ with $D(G(S)G(T))=D(G(S*T))\cap D(G(T))$, and $G(S)+G(T)\subseteq G(S+T)$.
\item[(iii)] $(\mathcal{G}(\psi)x$, $\mathcal{G}(-\psi^{\prime})x-\psi(0)Cx)\in G(-\delta')$.
\item[(iv)] If $\mathcal{G}$ is dense, then its generator is densely defined.
\end{itemize}
\end{prop}

The assertions (ii)-(vi) of \cite[Proposition 3.1.2]{knjigah} can be reformulated for pre-(C-UDS)'s of $\ast$-class in locally convex spaces.

\begin{prop}\label{kuki}
Let $\mathcal{G}$ be a pre-(C-UDS) of $\ast$-class. % $F:=E/\mathcal{N}(\mathcal{G})$
%and let $q$ be the corresponding canonical mapping $q:E\to F$.
Then the following holds:
\begin{itemize}
%\item[(i)] Let $H\in L(\mathcal{D}:L(F))$ be defined by $q\mathcal{G}(\varphi):=H(\varphi)q$
%for all $\varphi\in\mathcal{D}$ and let $\tilde{C}$ be a linear operator in $F$ defined by $\tilde{C}q:=qC$.
%Then $\tilde{C}\in L(F)$ and $\tilde{C}$ is injective.
%Moreover, $H$ is a $(\tilde{C}$-DS) in $F$.
\item[(i)] $C(\overline{\langle\mathcal{R}(\mathcal{G})\rangle})\subseteq\overline{\mathcal{R}(\mathcal{G})}$,
where $\langle\mathcal{R}(\mathcal{G})\rangle$
denotes the linear span of $\mathcal{R}(\mathcal{G})$.
\item[(ii)] Assume $\mathcal{G}$ is not dense and
$\overline{C\mathcal{R}(\mathcal{G})}=\overline{\mathcal{R}(\mathcal{G})}$.
Put $R:=\overline{\mathcal{R}(\mathcal{G})}$ and $H:=\mathcal{G}_{|R}$.
Then $H$ is a dense pre-($C_1$-UDS) of $\ast$-class on $R$ with $C_1=C_{|R}$.
\item[(iii)] The dual $\mathcal{G}(\cdot)^*$ is a pre-($C^*$-UDS) of $\ast$-class on $E^*$
and $\mathcal{N}(\mathcal{G}^*)=\overline{\mathcal{R}(\mathcal{G})}^{\circ}$.
\item[(iv)] If $E$ is reflexive,
then $\mathcal{N}(\mathcal{G})=\overline{\mathcal{R}(\mathcal{G}^*)}^{\circ}$.
\item[(v)]
The $\mathcal{G}^*$ is a ($C^*$-UDS) of $\ast$-class in $E^*$ iff $\mathcal{G}$ is a dense pre-(C-UDS) of $\ast$-class.
If $E$ is reflexive, then $\mathcal{G}^*$ is a dense pre-($C^*$-UDS) of $\ast$-class in $E^*$ iff $\mathcal{G}$ is a (C-UDS) of $\ast$-class.
\end{itemize}
\end{prop}

The following proposition has been recently proved in \cite{C-ultra} in the case that the operator $C$ is injective (cf. \cite[Proposition 2]{ki90}). By the proof of the statement in \cite{C-ultra}, it is clear that the injectivity of $C$ is superfluous.
\begin{prop}\label{kisinski}
Suppose that ${\mathcal G}\in {\mathcal D}^{\prime \ast}_{0}(L(E))$ and ${\mathcal G}(\varphi)C=C{\mathcal G}(\varphi),$ $\varphi \in {\mathcal D}^{\ast}$.
Then ${\mathcal G}$ is a pre-(C-UDS) of $\ast$-class iff
\begin{align*}
{\mathcal G}\bigl(\varphi^{\prime}\bigr){\mathcal G}(\psi)-{\mathcal G}(\varphi){\mathcal G}\bigl(\psi^{\prime}\bigr)=\psi(0){\mathcal G}(\varphi)C-\varphi(0){\mathcal G}
(\psi)C,\quad  \varphi,\ \psi \in {\mathcal D}^{\ast}.
\end{align*}
\end{prop}

In \cite{C-ultra}, we have recently proved that
every (C-UDS) of $\ast$-class in locally convex space is uniquely determined by its generator.
Contrary to the single-valued case, different pre-(C-UDS)'s of $\ast$-class
can have the same generator.

Next we give the following definition of an exponential pre-(C-UDS) of $\ast$-class.

\begin{defn}\label{qdfn}
Let ${\mathcal G}$ be a pre-(C-UDS) of $\ast$-class.
Then $\mathcal{G}$ is said to be an exponential
pre-(C-UDS) of $\ast$-class iff there exists $\omega\in\mathbb{R}$ such that $e^{-\omega t}\mathcal{G}
\in\mathcal{S}'^{\ast}(L(E))$.
We use the shorthand pre-(C-EUDS) of $\ast$-class to denote an exponential
pre-(C-UDS) of $\ast$-class.
\end{defn}

\begin{rem}\label{fundamentalna}
Suppose that ${\mathcal G}\in {\mathcal D}^{\prime \ast}_{0}(L(E))$, ${\mathcal G}(\varphi)C=C{\mathcal G}(\varphi),$ $\varphi \in {\mathcal D}^{\ast}$
and ${\mathcal A}$ is a closed MLO on $E$ satisfying that $\mathcal{G}(\varphi){\mathcal A}\subseteq {\mathcal A}{\mathcal G}(\varphi),$ $\varphi \in {\mathcal D}^{\ast}$ and
\begin{equation}\label{dkenk}
\mathcal{G}\bigl(-\varphi'\bigr)x-\varphi(0)Cx\in {\mathcal A}\mathcal{G}(\varphi)x,\quad x\in E,\ \varphi \in {\mathcal D}^{\ast}.
\end{equation}
\end{rem}
\begin{defn}
Let $\mathcal{G}$ be a pre-C-ultradistribution semigroup (pre-C-distribution semigroup). Then $\mathcal{G}$ is said to be a quasi-equicontnuous  exponential (short, (q)-exponential) pre-C-ultradistribution semigroup (pre-C-distribution semigroup) if for every $p\in\circledast$ and bounded subset $B\in E$ there exist $M_p\geq1$, ${\omega}_p\geq0$ and $q_p$ seminorm on $\mathcal{S}^{\ast}({\mathbb R})$ ($\mathcal{S}({\mathbb R})$) such that $$\sup\limits_{x\in B}p({\mathcal{G}}(\varphi)x)\leq M_pe^{{\omega}_p\cdot}q_p(\varphi),$$ for all $\varphi\in\mathcal{S}^{\ast}_0({\mathbb R})$ ($\varphi\in{\mathcal{S}}_0({\mathbb R})$). We use the shorthand pre-q-(C-EUDS) (pre-q-(C-EDS)).
\end{defn}

The following statements hold (see \cite{C-ultra}):
\begin{itemize}
\item[(i)] If ${\mathcal A}=A$ is single-valued, then ${\mathcal G}$ satisfies (C.S.1).
\item[(ii)] If ${\mathcal G}$ satisfies (C.S.2) holds, $C$ is injective and ${\mathcal A}=A$ is single-valued, then ${\mathcal G}$ is a (C-UDS) of $\ast$-class generated by $C^{-1}AC.$
\end{itemize}
As we have already seen, the conclusion from (ii) immediately implies that ${\mathcal A}=A$ must be single-valued and that
the operator $C$ must be injective.

Concerning the assertion (i), its validity is not true in multivalued case:
Let $C=I,$ let ${\mathcal A}\equiv E\times E,$ and let ${\mathcal G}\in {\mathcal D}^{\prime \ast}_{0}(L(E))$ be arbitrarily chosen. Then ${\mathcal G}$ commutes with ${\mathcal A}$ and
 (\ref{dkenk}) holds but ${\mathcal G}$ need not satisfy (C.S.1).

Concerning degenerate C-ultradistribution semigroups, exponential degenerate C-ultradistribution semigroups and degenerate (q-)exponential C-ultradistribution semigroups, we can give the following theorems (see \cite{quasirc}).
\begin{thm}\label{tempera-ultra}
\begin{itemize}
\item[(i)] Suppose that there exist $l>0,$ $\beta>0$ and $k>0,$ in the Beurling case,
resp., for every $l>0$ there exists $\beta_l>0$, in the Roumieu case, such that
$
\Omega^{(M_p)}_{l,\beta}:=\{\lambda\in\mathbb{C}:\Re\lambda\geq M(l|\lambda|)+\beta\}\subseteq\rho_{C}(A),
$ resp.
$
\Omega^{\{M_p\}}_{l,\beta_l}:=\{\lambda\in\mathbb{C}:\Re\lambda\geq M(l|\lambda|)+\beta_l\}\subseteq\rho_{C}(A),
$
the mapping $\lambda \rightarrow (\lambda-A)^{-1}Cx,$ $\lambda \in \Omega^{(M_p)}_{l,\beta},$ resp. $\lambda \in \Omega^{\{M_p\}}_{l,\beta_{l}},$ is continuous for every fixed element $x\in E,$
and the operator family
$
\{e^{-M(kl|\lambda|)}(\lambda-A)^{-1}C : \lambda\in\Omega^{(M_p)}_{l,\beta}\}\subseteq L(E),
$ resp. $
\{e^{-M(l|\lambda|)}(\lambda-A)^{-1}C : \lambda\in\Omega^{\{M_p\}}_{l,\beta_{l}}\}\subseteq L(E),
$
is equicontinuous. Denote by $\Gamma ,$ resp. $\Gamma_{l},$ the upwards oriented boundary of $\Omega^{(M_p)}_{l,\beta},$ resp. $\Omega^{\{M_p\}}_{l,\beta_{l}}.$ Define, for every $x\in E$ and $\varphi\in\mathcal{D}^{\ast},$ the element
$
{\mathcal G}(\varphi)x$ with
\begin{align} \mathcal{G}(\varphi)x:=(-i)\int_{\Gamma}\hat{\varphi}(\lambda)(\lambda-A)^{-1}Cx\,d\lambda,
\;\;x\in E,\;\varphi\in\mathcal{D},
\end{align} in the Beurling case; in the Roumieu case, for every number $k>0$ and for every function $\varphi \in {\mathcal D}^{\{M_{p}\}}_{[-k,k]},$
we define the element ${\mathcal G}(\varphi)x$ in the same way as above, with the contour $\Gamma$ replaced by $\Gamma_{l(k)}$.
Then ${\mathcal G}\in {\mathcal D}^{\prime \ast}_{0}(L(E))$ is boundedly equicontinuous,
${\mathcal G}(\varphi)C=C{\mathcal G}(\varphi),$ $\varphi \in {\mathcal D}^{\ast},$
$\mathcal{G}(\varphi)A\subseteq A{\mathcal G}(\varphi),$ $\varphi \in {\mathcal D}^{\ast}$ and $A\mathcal{G}(\varphi)x=\mathcal{G}\bigl(-\varphi'\bigr)x-\varphi(0)Cx,\quad x\in E,\ \varphi \in {\mathcal D} \ \ (\varphi \in {\mathcal D}^{\ast})$. Then, ${\mathcal G}$ is a pre-(C-UDS) of $\ast$-class.
\item[(ii)]
Suppose that $A$ is a closed linear operator on $E$ satisfying that there exist $a\geq 0$ such that $
\{\lambda\in\mathbb{C}:\Re\lambda>a\}\subseteq\rho_{C}(A)$ and the mapping $\lambda \mapsto (\lambda-A)^{-1}Cx,$ $\Re \lambda>a$ is continuous for every fixed element $x\in E.$
Suppose that there exists a number $k>0,$ in the Beurling case, resp., for every number $k>0,$
in the Roumieu case,
such that the operator family
$\{e^{-M(k|\lambda|)}(\lambda-A)^{-1}C : \Re \lambda>a\}\subseteq L(E)$ is equicontinuous. Set
$$
{\mathcal G}(\varphi)x=(-i)\int\limits_{\bar{a}-i\infty}^{\bar{a}
+i\infty}\hat{\varphi}(\lambda)\bigl(\lambda-A\bigr)^{-1}Cx\,d\lambda,\quad x\in E,\ \varphi\in\mathcal{D}^{\ast}.
$$
Then ${\mathcal G}\in {\mathcal D}^{\prime \ast}_{0}(L(E))$ is boundedly equicontinuous, $e^{-\omega t}G\in\mathcal{S}^{\prime \ast}(L(E))$ for all $\omega >a,$
${\mathcal G}(\varphi)C=C{\mathcal G}(\varphi),$ $\varphi \in {\mathcal D}^{\ast},$
$\mathcal{G}(\varphi)A\subseteq A{\mathcal G}(\varphi),$ $\varphi \in {\mathcal D}^{\ast}$ and $A\mathcal{G}(\varphi)x=\mathcal{G}\bigl(-\varphi'\bigr)x-\varphi(0)Cx,\quad x\in E,\ \varphi \in {\mathcal D} \ \ (\varphi \in {\mathcal D}^{\ast})$. Then, ${\mathcal G}$ is a pre-(C-EUDS) of $\ast$-class.
\end{itemize}
\end{thm}

\begin{rem} By J. Chazarain \cite{cha}, we  define $(M_p)$-ultralogarithmic region ${\Lambda}_{\alpha,\beta,l}$ of type $l$ as $${\Lambda}_{\alpha,\beta,l}=\{\lambda\in{\mathbb C}\,\, :\,\, \Re\lambda\geq\alpha M(l|\Im\lambda|)+\beta\},$$ for $\alpha,\beta>0$, $l\in{\mathbb R}$. The first part of the Theorem \ref{tempera-ultra} can be reformulated with the region ${\Omega}_{l,\beta}^{(M_p)}$ replaced by ${\Lambda}_{\alpha,\beta,l}$.\end{rem}

Let $\bar{\alpha}>\alpha.$ By $\Gamma_{l}$ ($\Gamma_{\bar{\alpha}}$) we donte the upwards oriented boundary of the ultra-logarithmic region $\Lambda_{\alpha,\beta, l}$ (the right line connecting the points $\bar{\alpha}-i\infty$ and $\bar{\alpha}+i\infty$) and let
\begin{equation}\label{polugrupa}
{\mathcal G}(\varphi)x:=(-i)\int_{\Gamma_{l}\ (\Gamma_{\bar{\alpha}})}\hat{\varphi}(\lambda)(\lambda-A)^{-1}x\, d\lambda,\quad x\in E,\quad \varphi \in {\mathcal D}^{(M_{p})}.
\end{equation}
The abstract Beurling space of $(M_p)$ class
associated to a closed linear operator $A$ is defined as in \cite{ci1}.
Following \cite{ci1}, we put $E^{(M_p)}(A):=$projlim$_{h\to+\infty}E^{(M_p)}_h(A)$, where
$$
E^{(M_p)}_h(A):=\Biggl\{x\in D_{\infty}(A):\|x\|^{(M_p)}_{h,q}=\sup_{p\in\mathbb{N}_0}\frac{h^{p}q\bigl(A^px\bigr)}{M_p}<\infty\mbox{ for all }h>0\mbox{ and } q\in \circledast\Biggr\}.
$$
Then, for each number $h>0$ the calibration $(\|\cdot\|^{(M_p)}_{h,q})_{q\in \circledast}$ induces a Hausdorff sequentially complete locally convex space on $E^{(M_p)}_h(A),$
$E^{(M_p)}_{h'}(A)\subseteq E^{(M_p)}_h(A)$
provided $0<h<h'<\infty ,$  and the spaces  $E^{(M_p)}_h(A)$  and  $E^{(M_p)}(A)$ are continuously  embedded in $E$.

\begin{thm}\label{grace}
Let $A$ be a closed linear operator $A$ and there exist constants $l\geq 1$, $\alpha>0,$ $\beta>0$ and $k>0$
such that $
\Lambda_{\alpha,\beta, l}\subseteq\rho (A)$ ($RHP_{\alpha}\equiv \{\lambda \in {\mathbb C} : \Re \lambda>\alpha\}\subseteq \rho (A)$). Let for each seminorm $q\in \circledast$ there exist a number $c_{q}>0$ and a seminorm $r\in \circledast$ such that
\begin{equation}\label{but-surf}
q\Bigl( \bigl(\lambda-A\bigr)^{-1}x \Bigr)\leq c_{q}e^{M(kl|\lambda|)}r(x),\quad x\in E,\ \lambda \in \Lambda_{\alpha,\beta, l} \ \ \bigl(RHP_{\alpha} \bigr).
\end{equation} Moreover, assume that ${\mathcal G},$ defined through \emph{(\ref{polugrupa})}, is a (UDS) ((EUDS)) of Beurling class generated by $A$ (i.e., that ${\mathcal G}$ satisfies \emph{(C.S.2)}), and that  $(M_{p})$ satisfies
\emph{(M.1)}, \emph{(M.2)} and \emph{(M.3)}.
Then the abstract Cauchy problem $(ACP)$
has a unique solution $u(t)$ for all $x\in E^{(M_p)}(A)$.
\end{thm}
\begin{rem}
By the discussion made before the last two theorems, we can say that $A$ can not be a multivalued linear operator in Theorem \ref{tempera-ultra} and Theorem \ref{grace}.
\end{rem}

%Concerning the assertion (iii) in multivalued case, we can prove that the admissibility of state space $E$ implies that for each
%$x\in {\mathcal N}({\mathcal G})$
%there exist an integer $k\in {\mathbb N}$ and a finite sequence $(y_{i})_{0\leq i \leq k-1}$ in $D({\mathcal A})$ such that
%$y_{i}\in {\mathcal A}y_{i+1}$ ($0\leq i \leq k-1$) and $Cx\in {\mathcal A}y_{0}\subseteq {\mathcal A}^{k+2}0.$

Now we will reconsider some conditions (originally introduced by J. L. Lions \cite{li121}, for the definition of dense distribution semigroups and for ultradistribution case the conditions in \cite{C-ultra}) in our new framework.
Suppose that $\mathcal{G}\in\mathcal{D}'^{\ast}_0(L(E))$ and ${\mathcal G}$ commutes with $C$.
Like in the case of degenerate $C$-distribution semigroups (see \cite{degdis}), we analyze the following conditions for ${\mathcal G}$:
\begin{itemize}
\item[$(d_1)$] $\mathcal{G}(\varphi*\psi)C=\mathcal{G}(\varphi)\mathcal{G}(\psi)$, $\varphi,\,\psi\in\mathcal{D}^{\ast}_0$,
\item[$(d_3)$] $\mathcal{R}(\mathcal{G})$ is dense in $E$,
\item[$(d_4)$] for every $x\in\mathcal{R}(\mathcal{G})$, there exists a function $u_x\in C([0,\infty):E)$ so that
$u_x(0)=Cx$ and $\mathcal{G}(\varphi)x=\int_0^{\infty}\varphi(t)u_x(t)\,dt$, $\varphi\in\mathcal{D}^{\ast}$,
\item[$(d_5)$] $(Cx,\mathcal{G}(\psi)x)\in G(\psi_+)$, $\psi\in\mathcal{D}^{\ast},$ $x\in E$.
\end{itemize}
We will discuss the connections of the previously given conditions, $(d_1)$, $(d_2)$, $(d_3)$, $(d_4)$ and $(d_5)$.
Let $\mathcal{G}\in\mathcal{D}'^{\ast}_0(L(E))$ be a pre-(C-UDS) of $\ast$-class. Then
$\mathcal{G}$
satisfies $(d_1)$ and from previously $\mathcal{G}$
satisfies $(d_5)$. Also, by the proof of \cite[Proposition 3.1.24]{knjigah}, we have that $\mathcal{G}$
also satisfies $(d_4).$ On the other hand, it is well known that $(d_1),$ $(d_4)$ and (C.S.2) taken together do not imply
(C.S.1), even in the case that $C=I;$ see e.g. \cite[Remark 3.1.20]{knjigah}. Furthermore, if $(d_1),$ $(d_3)$ and $(d_4)$ hold then  $(d_5)$ holds, as well. To prove this, fix $x\in\mathcal{R}(\mathcal{G})$ and $\varphi\in\mathcal{D}^{\ast}$. Then it suffices to show that
$(Cx,\mathcal{G}(\varphi)x)\in G(\varphi_+)$.  Suppose that $(\rho_n)$ is a regularizing sequence
and $u_x(t)$ is a function appearing in the formulation of the property $(d_4)$.
From the proof of \cite[Proposition 3.1.19]{knjigah}, for every $\eta\in\mathcal{D}_{0}^{\ast}$, we have
\begin{align*}
\mathcal{G}(\rho_n)\mathcal{G}(\varphi_+*\eta)x&=\mathcal{G}((\varphi_+*\rho_n)*\eta)Cx
=\mathcal{G}(\eta)\mathcal{G}(\varphi_+*\rho_n)x\\
&=\mathcal{G}(\eta)\int\limits_0^{\infty}(\varphi_+*\rho_n)(t)u_x(t)\,dt
\\
& \to \mathcal{G}(\eta)\int_0^{\infty}\varphi(t)u_x(t)\,dt=\mathcal{G}(\eta)\mathcal{G}(\varphi)x,\;n\to\infty;\\
\mathcal{G}(\rho_n)\mathcal{G}(\varphi_+*\eta)x
&=\mathcal{G}(\varphi_+*\eta*\rho_n)Cx\to\mathcal{G}(\varphi_+*\eta)Cx,\;n\to\infty.
\end{align*}
Hence, $\mathcal{G}(\varphi_+*\eta)Cx=\mathcal{G}(\eta)\mathcal{G}(\varphi)x$ and $(d_5)$ holds, as claimed. On the other hand,
$(d_1)$ is a very simple consequence of $(d_5).$ To see this, observe that for each $\varphi \in\mathcal{D}_0^{\ast}$ and $\psi \in\mathcal{D}^{\ast}$ we have $\psi_{+}*\varphi=\psi \ast_0 \varphi =\varphi \ast_0 \psi,$ so that $(d_5)$ is equivalent to say that
$\mathcal{G}(\varphi \ast_0 \psi)C=\mathcal{G}(\varphi)\mathcal{G}(\psi)$ ($\varphi \in\mathcal{D}_0^{\ast},$ $\psi \in\mathcal{D}^{\ast}$). In particular,
\begin{align}\label{line988}
\mathcal{G}(\varphi)\mathcal{G}(\psi)=\mathcal{G}(\psi)\mathcal{G}(\varphi),\quad \varphi \in\mathcal{D}_0^{\ast},\ \psi \in\mathcal{D}^{\ast}.
\end{align}
Now, let $(d_5)$ holds, $\varphi \in\mathcal{D}_0^{\ast}$ and $\psi,\ \eta \in\mathcal{D}^{\ast}.$ Note that $\psi_{+}\ast \eta_{+} \ast \varphi=(\psi \ast_{0}\eta)_{+} \ast \varphi$. Then (cf. also \cite[Remark 3.13]{ku112}):
\begin{align}
\notag \mathcal{G}(\varphi)\mathcal{G}(\eta)\mathcal{G}(\psi)&= C\mathcal{G}(\eta_{+} \ast \varphi)\mathcal{G}(\psi)
\\\notag &=C\mathcal{G}(\psi_{+}\ast \eta_{+} \ast \varphi)=C\mathcal{G}\bigl( (\psi \ast_{0}\eta)_{+} \ast \varphi \bigr)C
\\\label{line98888} &=C\mathcal{G}(\varphi)\mathcal{G}( \psi \ast_{0}\eta )=\mathcal{G}(\varphi)\mathcal{G}( \psi \ast_{0}\eta )C.
\end{align}
By (\ref{line988})-(\ref{line98888}), we get
\begin{align}\label{line988888}
\mathcal{G}(\eta)\mathcal{G}(\psi)\mathcal{G}(\varphi)=\mathcal{G}( \psi \ast_{0}\eta )C\mathcal{G}(\varphi).
\end{align}
By (\ref{line988})-(\ref{line988888}), we have the following conclusions:
\begin{itemize}
\item[(i)] $(d_5)$ and $(d_3)$ together imply (C.S.1); in particular, $(d_1),$ $(d_3)$ and $(d_4)$ together imply  (C.S.1). This is an extension of \cite[Proposition 3.1.19]{knjigah}.
\item[(ii)] $(d_5)$ and $(d_2)$ together imply that ${\mathcal G}$ is a (C-UDS) of $\ast$-class; in particular, ${\mathcal A}=A$ must be single-valued and $C$ must be injective.
\end{itemize}
On the other hand, $(d_5)$ does not imply (C.S.1) even in the case that $C=I.$ A simple counterexample is $\mathcal{G}\in\mathcal{D}'^{\ast}_0(L(E))$
given by  ${\mathcal G}(\varphi)x:=\varphi(0)x,$ $x\in E,$ $\varphi \in {\mathcal D}^{\ast}$.

The exponential region $E(a,b)$ has been defined for the first time by W. Arendt, O. El--Mennaoui and V. Keyantuo in \cite{a22}:
$$
E(a,b):=\Bigl\{\lambda\in\mathbb{C}:\Re\lambda\geq b,\:|\Im\lambda|\leq e^{a\Re\lambda}\Bigr\} \ \ (a,\ b>0).
$$

\begin{rem}\label{filipa-1}
%\begin{itemize}
%\item[(i)]
%If $C$ is injective, ${\mathcal A}=A$ is single-valued, $\rho_{C}(A) \subseteq E(a,b)$ and $F(\lambda)=(\lambda -{\mathcal A})^{-1}C,$ $\lambda \in E(a,b),$ then ${\mathcal G}$ is a (C-DS) generated by $C^{-1}AC$
%(\cite{C-ultra}).
%Even in the case that $C=I,$ the integral generator ${\mathcal A}$ of ${\mathcal G},$ in multivalued case, can strictly contain $C^{-1}{\mathcal A}C;$ see Remark \ref{xcvbnm-prim}(i).
%\item[(ii)] Let ${\mathcal A}$ be a closed MLO, let $C$ be injective and commute with ${\mathcal A},$ and let $\rho_{C}({\mathcal A}) \subseteq E(a,b).$ Then the choice $F(\lambda)=(\lambda -{\mathcal A})^{-1}C,$ $\lambda \in E(a,b)$
%is always possible; in this case, we have ${\mathcal A}0\subseteq
%N({\mathcal G}(\varphi)),$ $\varphi \in {\mathcal D}$  (\cite{FKP}).
Suppose that there exist $l>0,$ $\beta>0$ and $k>0,$ in the Beurling case,
resp., for every $l>0$ there exists $\beta_l>0$, in the Roumieu case, such that the assumptions of \cite[Theorem 4.15]{degdis} hold with the exponential region $E(a,b)$ replaced with the region
$
\Omega^{(M_p)}_{l,\beta}:=\{\lambda\in\mathbb{C}:\Re\lambda\geq M(l|\lambda|)+\beta\},
$ resp.
$
\Omega^{\{M_p\}}_{l,\beta_l}:=\{\lambda\in\mathbb{C}:\Re\lambda\geq M(l|\lambda|)+\beta_l\}.
$
Define $
{\mathcal G}
$ similarly as above. Then $\mathcal{G}\in {\mathcal D}^{\prime \ast}_{0}(L(E)),$ $\mathcal{G}$ commutes with $C$ and ${\mathcal A},$ and
(\ref{dkenk}) holds. But, in the present situation, we do not know whether
$\mathcal{G}$ has to satisfy (C.S.1) in degenerate case. This is an open problem we would like to address to our readers.
%\end{itemize}
\end{rem}

%Local integrated semigroups generated by multivalued linear operators (see e.g. \cite[Example 3.2.11(i)]{FKP}) can be used for construction of pre-(DS)'s.
Let we mention that in \cite[Theorem 3.2.21]{FKP} and \cite[Example 3.2.23]{FKP}, are investigated the entire solutions
of backward heat Poisson equation and is showed the existence of
an entire $C$-regularized semigroup ($C\in L(L^{p}(\Omega))$ non-injective) generated by the multivalued linear operator $\Delta \cdot m(x)^{-1}$ in $L^{p}(\Omega),$
where $\Omega$ is a bounded domain in ${\mathbb R}^{n}$.
%This example can serve us to construct an important example of a pre-(C-DS); cf. also \cite[Example 3.24]{catania}. Examples of exponentially bounded integrated semigroups generated by multivalued linear operators can be found in \cite[Chapter II-III, Section 5.8]{faviniyagi} and these examples can be used for construction of exponential pre-(DS)'s.

An example of exponential degenerate ultradistribution semigroup of Beurling class can be given by using the consideration  from \cite[Example 3.25]{catania}. By Proposition \ref{kuki}(iii), the duals of
non-dense (C-UDS)'s of $\ast$-class serve as examples of
pre-(C$^{\ast}$-UDS)'s of $\ast$-class, as well.

\end{document}